\documentstyle[twoside]{article}
\title{\bf Asymptotics of a sum of modified Bessel functions with a non-linear argument}
\author{\sc R. B. Paris\footnote{E-mail address:\ \ {\tt r.paris@abertay.ac.uk}}\\
\\
{\em Division of Computing and Mathematics,}\\
{\em Abertay University, Dundee DD1 1HG, UK}\\
}

\begin{document}
\newcommand{\bee}{\begin{equation}}
\newcommand{\ee}{\end{equation}}
\def\f#1#2{\mbox{${\textstyle \frac{#1}{#2}}$}}
\def\dfrac#1#2{\displaystyle{\frac{#1}{#2}}}
\newcommand{\fr}{\frac{1}{2}}
\newcommand{\fs}{\f{1}{2}}
\newcommand{\g}{\Gamma}
\newcommand{\om}{\omega}
\newcommand{\br}{\biggr}
\newcommand{\bl}{\biggl}
\newcommand{\ra}{\rightarrow}
\renewcommand{\topfraction}{0.9}
\renewcommand{\bottomfraction}{0.9}
\renewcommand{\textfraction}{0.05}
\newcommand{\mcol}{\multicolumn}
\date{}
\maketitle
\pagestyle{myheadings}
\markboth{\hfill {\it R.B. Paris} \hfill}
{\hfill {\it A modified Bessel function sum } \hfill}
\begin{abstract} 
We examine the sum of modified Bessel functions with argument depending non-linearly on the summation index given by
\[S_{\nu,p}(a)=\sum_{n\geq 1} (\fs an^p)^{-\nu} K_\nu(an^p)\qquad (a>0,\ 0\leq\nu<1)\]
as the parameter $a\to 0+$, where $p$ denotes an integer satisfying $p\geq 2$. This extends previous work for the cases $p=1$ (linear) and $p=2$ (quadratic).
The expansion as $a\to0+$ consists of an infinite number of asymptotic sums involving the Riemann zeta function, which when optimally truncated lead to remainder terms that are exponentially small in the parameter $a$. The number of these exponentially small terms associated with each optimally truncated asymptotic sum is found to increase with $p$. 
\vspace{0.4cm}

\vspace{0.3cm}

\end{abstract}

\vspace{0.2cm}

\noindent $\,$\hrulefill $\,$

\vspace{0.2cm}

\begin{center}
{\bf 1. \  Introduction}
\end{center}
\setcounter{section}{1}
\setcounter{equation}{0}
\renewcommand{\theequation}{\arabic{section}.\arabic{equation}}
We consider the asymptotic expansion of the sum
\bee\label{e11}
S_{\nu,p}(a)=\sum_{n\geq 1} (\fs an^p)^{-\nu} K_\nu(an^p),\qquad a>0
\ee
as the parameter $a\to 0+$, where $p\geq 2$ is an integer. The function $K_\nu(z)$ is the modified Bessel function of the second kind of order $\nu$ that we restrict to satisfy $0\leq\nu<1$. The argument of the Bessel function depends non-linearly on the summation index $n$, thereby extending the sum considered in \cite{P18} with argument depending quadratically on $n$. The sum converges without restriction on $\nu$ on account of the exponential decay of the modified Bessel function for large argument.

When $\nu=\fs$ the $K$-Bessel function reduces to a simple exponential and the above sum is then expressible in terms of the generalised Euler-Jacobi series 
\bee\label{e12}
S_{\frac{1}{2},p}(a)=\frac{\sqrt{\pi}}{2}\bl(\frac{2}{a}\br)^{\!\nu+\frac{1}{2}} \sum_{n\geq 1}\frac{e^{-an^p}}{n^w},\qquad w=p(\nu+\fs).
\ee
Sums of this type have a long history dating back to Poisson in the classical case $p=2$, $w=0$; see \cite[p.~124]{WW}. They were studied by Ramanujan  who showed that the expansion as $a\to 0$ consisted of an asymptotic sum involving the Riemann zeta function; see \cite[Chapter 1]{Berndt}. Detailed studies of the Euler-Jacobi series in the case $w=0$ have been carried out in \cite{Kow} and \cite[Section 8.1]{PK}, and more recently in \cite{P16}, where it was established that as $a\to0$ an increasing number of exponentially small contributions arise as the parameter $p$ increases.

In the limit $a\to 0$ the resulting slow convergence of the terms in the sum $S_{\nu,p}(a)$ renders computation more difficult. It will be shown that the asymptotic expansion of $S_{\nu,p}(a)$ in this limit consists of a finite series in ascending powers of $a$ with coefficients involving the Riemann zeta function, together with an infinite number of asymptotic sums associated with increasingly subdominant exponentially small terms. 

\vspace{0.6cm}

\begin{center}
{\bf 2. \ An expansion for $S_{\nu,p}(a)$}
\end{center}
\setcounter{section}{2}
\setcounter{equation}{0}
\renewcommand{\theequation}{\arabic{section}.\arabic{equation}}
We start with the Mellin-Barnes integral representation\footnote{There is an error in the sector of validity in \cite[(10.32.13)]{DLMF} and also in \cite[p.~114]{PK}. The quantity $c$ will be used as a generic parameter that can vary according to each integral.} \cite[(10.32.13)]{DLMF}
\[(\fs x)^{-\nu} K_\nu(x)=\frac{1}{4\pi i}\int_{c-\infty i}^{c+\infty i}\g(s) \g(s-\nu) (\fs x)^{-2s} ds\qquad (|\arg\,x|<\fs\pi),\]
where $c>\max\{0,\nu\}$ so that the integration path lies to the right of all the poles of the integrand.  Then it follows that
\bee\label{e21}
{\cal S}_{\nu,p}(a)=\sum_{n\geq 1}(\fs an^p)^{-\nu} K_\nu(an^p)=\frac{1}{4\pi i}\int_{c-\infty i}^{c+\infty i} \g(s)\g(s-\nu)(\fs a)^{-2s} \zeta(2ps)\,ds
\ee
for $|\arg\,a|<\fs\pi$, where $\zeta(s)$ is the Riemann zeta function and $c>\max\{\nu, 1/(2p)\}$.

For simplicity in presentation we shall suppose that $0\leq\nu<1$; the case $\nu\geq 1$, which involves additional poles in $\Re (s)>0$, is straightforward and can be handled as in \cite{P18}.
The integrand in (\ref{e21}) has a pole at $s=1/(2p)$ resulting from $\zeta(2ps)$. The integrand also has poles at $s=\nu-n$, $n=0, 1, 2, \ldots$ together with a pole at $s=0$, the other poles of $\g(s)$ at $s=-1, -2, \ldots$ being cancelled by the trivial zeros of $\zeta(2ps)$ at $ps=-1, -2, \ldots\, $. When $\nu=1/(2p)$ the pole at $s=\nu$ is double and when $\nu=0$ the pole at $s=0$ is double.

We consider the integral taken round the rectangular contour with vertices at $c\pm iT$, $-d\pm iT$, where $0<d<1-\nu$.
The contribution from the upper and lower sides $s=\sigma\pm iT$, $-d\leq \sigma\leq c$ as $T\to\infty$ can be estimated by use of the standard results
\[|\g(\sigma\pm it)|\sim \sqrt{2\pi}t^{\sigma-\frac{1}{2}}e^{-\frac{1}{2}\pi t}\qquad (t\to+\infty),\]
which follows from Stirling's formula for the gamma function, and \cite[p.~25]{Iv}
\[|\zeta(\sigma\pm it)|=O(t^{{\hat\mu}(\sigma)}\log^\beta t)\qquad (t\to+\infty)\]
where ${\hat\mu}(\sigma)=0$ ($\sigma>1$), $\fs-\fs\sigma$ ($0\leq\sigma\leq1$), $\fs-\sigma$ ($\sigma\leq0$) and $\beta=0$ ($\sigma>1$), 1 ($\sigma\leq 1)$. Then it follows that
\[\g(\sigma\pm it) \g(\sigma-\nu\pm it)(\fs a)^{-2\sigma\mp 2it}=O(t^{2\sigma-\nu-1} e^{-\pi t})\]
as $t\to +\infty$.
Hence the modulus of the integrand on these horizontal paths is $O(T^{\xi} \log\,T e^{-\Delta T})$ as $T\to\infty$, where $\xi=\sigma+{\hat\mu}(2p\sigma)-\nu-1$ and $\Delta=\pi-2|\arg\,a|$. Taking into account the different forms of ${\hat\mu}(2p\sigma)$, we see that the contribution from these paths vanishes as $T\to\infty$ provided $|\arg\,a|<\fs\pi$.

We first displace the integration path to the left over the poles in $\Re (s)\geq 0$. We make use of the results that $\zeta(s)$ has residue equal to 1 at $s=1$, $\zeta(0)=-\fs$, $\zeta'(0)=-\fs\log\,2\pi$ and
\[\g(\alpha+\epsilon)=\g(\alpha)\{1+\epsilon \psi(\alpha)+O(\epsilon^2)\},\qquad \zeta(1+\epsilon)=\epsilon^{-1}\{1+\epsilon\gamma+O(\epsilon^2)\} \qquad (\epsilon\to 0),\]
where $\psi(\alpha)$ denotes the logarithmic derivative of the gamma function and $\gamma=0.577215\ldots$ is the Euler-Mascheroni constant. Then we find after routine calculations the residue contribution in $\Re (s)\geq 0$ when $0\leq\nu<1$ given by
\bee\label{e2h}
{\cal H}_{\nu,p}(a)=\left\{\begin{array}{ll} \displaystyle{\frac{1}{2}\g(\nu)(\fs a)^{-2\nu} \zeta(2p\nu)+\frac{1}{4p}\g\bl(\frac{1}{2p}\br)\g\bl(\frac{1}{2p}-\nu\br)(\fs a)^{-1/p}+\frac{\pi\mbox{cosec}\,\pi\nu }{4\g(1+\nu)}} & (\nu\neq 0, \nu_*)\\
\\
\displaystyle{\frac{1}{4p}\g^2\bl(\frac{1}{2p}\br) (\fs a)^{-1/p}+\frac{1}{2}\{\gamma+\log\,\fs a-p \log\,2\pi\}} & (\nu=0)\\
\\
\displaystyle{\frac{1}{4p}\g\bl(\frac{1}{2p}\br) (\fs a)^{-1/p}\{\psi(\f{1}{2p})-2\log\,\fs a+(2p-1)\gamma\}+\!\frac{\pi\mbox{cosec}\,\pi/(2p)}{4\g(1+\frac{1}{2p})}}
& (\nu=\nu_*),
\end{array}\right.\ee
where $\nu_*=1/(2p)$. Then we have
\bee\label{e23}
{\cal S}_{\nu,p}(a)-{\cal H}_{\nu,p}(a)=\frac{1}{4\pi i}\int_{-c-\infty i}^{-c+\infty i}\g(s) \g(s-\nu) (\fs a)^{-2s} \zeta(2ps)\,ds,
\ee
where $0<c<1-\nu$. 

We now employ the functional relation for $\zeta(s)$ in the form \cite[(25.4.2)]{DLMF}
\[\zeta(s)=2^s\pi^{s-1} \g(1-s)\zeta(1-s) \sin \fs\pi s\]
and make the change of variable $s\to -s$ to obtain
\[{\cal S}_{\nu,p}(a)-{\cal H}_{\nu,p}(a)=-\frac{\pi}{4\pi i}\int_{c-\infty i}^{c+\infty i} \frac{\g(1+2ps)\zeta(1+2ps)}{\g(1\!+\! s)\g(1\!+\!s\!+\!\nu)}\,\bl(\frac{a}{2(2\pi)^p}\br)^{\!\!2s}\,\frac{\sin \pi p s}{\sin \pi s\,\sin \pi(s\!+\!\nu)}\,ds\]
\[=\sum_{k\geq 1}\frac{1}{k}\,J_k(a),\]
where
\bee\label{e24}
J_k(a)=-\frac{\pi}{4\pi i}\int_{c-\infty i}^{c+\infty i}
\frac{\g(1\!+\!2ps)}{\g(1\!+\! s)\g(1\!+\!s\!+\!\nu)}\,
\bl(\frac{a}{2(2\pi k)^p}\br)^{\!\!2s}\,\frac{\sin \pi ps}{\sin \pi s\,\sin\pi(s\!+\!\nu)}\,ds
\ee
with $0<c<1-\nu$.
Here we have used the series expansion for $\zeta(1+2ps)$ since $\Re (s)>0$ on the integration path
and reversal of the order of summation and integration is justified by absolute convergence.

A similar estimate of the integrand on horizontal paths $\Im (s)=\pm T$ as employed in the residue evaluation in (\ref{e2h}) shows that the integration path in (\ref{e24}) can be displaced to the right
over the simple poles at $s=n-\nu$, $1\leq n\leq N_k-1$, where $N_k$ is (for the moment) an arbitrary positive integer. Then 
\bee\label{e25}
J_k(a)=\frac{\sin \pi p\nu}{2\sin \pi\nu} \sum_{n=1}^{N_k-1}(-1)^{pn} \frac{(2p(n-\nu))!}{n! (n-\nu)!}\,\bl(\frac{a}{2(2\pi k)^p}\br)^{\!2n-2\nu}
+R_k(a;N_k),
\ee
where the remainder is given by
\bee\label{e26}
R_k(a;N_k)=-\frac{\pi}{4\pi i}\int_{{\cal L}_N}
\frac{\g(1\!+\!2ps)}{\g(1\!+\!s) \g(1\!+\!s\!+\!\nu)}\,\bl(\frac{a}{2(2\pi k)^p}\br)^{\!2s}\, \frac{\sin \pi ps}{\sin \pi s\,\sin \pi(s\!+\!\nu)}\,ds
\ee
and ${\cal L}_N$ denotes the rectilinear path $(-c\!+\!N_k\!-\nu\!-\!\infty i, -c\!+\!N_k\!-\!\nu\!+\!\infty i)$ with $0<c<1$.
 
We therefore have the following result.
\newtheorem{theorem}{Theorem}
\begin{theorem}$\!\!\!.$\  Let $0\leq\nu<1$ and the sum $S_{\nu,p}(a)$ be as defined in (\ref{e11}). Then, for integer $p\geq 2$,
\bee\label{e27}
{\cal S}_{\nu,p}(a)={\cal H}_{\nu,p}(a)+\sum_{k\geq 1}\frac{1}{k}\bl\{\frac{\sin \pi p\nu}{2\sin \pi\nu}\sum_{n=1}^{N_k-1} (-)^{pn}\frac{(2p(n-\nu))!}{n! (n-\nu)!}\,\bl(\frac{a}{2(2\pi k)^p}\br)^{\!2n-2\nu}\!\!
+R_k(a;N_k)\br\}
\ee
for $|\arg\,a|<\fs\pi$, where ${\cal H}_{\nu,p}(a)$  and the remainder term $R_k(a;N_k)$ are defined in (\ref{e2h}) and (\ref{e26}).
\end{theorem}
The contribution to ${\cal S}_{\nu,p}(a)$ from $\Re (s)<0$ in the integral in (\ref{e23}) has been decomposed into a $k$-sequence of component asymptotic series with scale $2(2\pi k)^p/a$, each associated with its own truncation index $N_k$ and remainder term $R_k(a;N_k)$. When $p=2$ this reduces to the expression obtained in (2.7) of \cite{P18}. 

The indices $N_k$ ($k=1, 2, \ldots$) are, for the moment, unspecified but will be subsequently chosen to be their optimal truncation values. In the remainder of this paper our attention will be devoted to an analysis of the remainder term $R_k(a;N_k)$, where we restrict the parameter $a$ to be positive and consider the limit $a\to0+$.
\vspace{0.6cm}

\begin{center}
{\bf 3. \ Expansion of the remainder term $R_k(a;N_k)$ for $a\to0+$}
\end{center}
\setcounter{section}{3}
\setcounter{equation}{0}
\renewcommand{\theequation}{\arabic{section}.\arabic{equation}}
Let us define the quantities
\bee\label{e30}
\kappa=2(p-1),\quad h=(2p)^{2p},\quad\vartheta=-\nu,\quad X_k=\kappa \bl(\frac{2}{a}\,\bl(\frac{\pi k}{p}\br)^p\br)^{1/(p-1)},
\ee
where $X_k>0$.
The indices $N_k$ will now be chosen to correspond to the optimal truncation values (that is, truncation just before the numerically smallest term) of the $k$-component asymptotic series in (\ref{e25}). On the assumption that $N_k\gg 1$, these are determined by the requirement
\[\bl(\frac{2}{a}(2\pi k)^p\br)^{\!2} \simeq \frac{(2p(N_k+1)-2p\nu)!}{(2pN_k-2p\nu)!}\,\frac{1}{(N_k+1)(N_k+1-\nu)}\] 
\[\hspace{0.5cm}=\frac{(2pN_k)^{2p}}{N_k^2}\,\bl\{1+\frac{{\cal A}}{N_k}+O(N_k^{-2})\br\},\]
where
\bee\label{e300}
{\cal A}=\frac{1}{2p}\sum_{r=1}^{2p}(r-2p\nu)+\nu-2=p-\f{3}{2}-(\kappa+1)\nu.
\ee
Then we find
\[\bl(\frac{2}{a}\,\bl(\frac{\pi k}{p}\br)^{\!\!p}\br)^{\!\!1/(p-1)}=N_k+\kappa^{-1}{\cal A}+O(N_k^{-1}),\]
so that the optimal truncation indices are specified by
\bee\label{e32}
X_k=\kappa N_k+{\cal A}+\alpha_k,
\ee
where $|\alpha_k|$ is bounded. It is seen that as $a\to 0+$ the truncation indices $N_k\to+\infty$, ($k=1,2, \ldots$).
 
The variable $s$ in the quotient of gamma functions in (\ref{e26}) is uniformly large on the displaced integration path ${\cal L}_N$ since $N_k\to+\infty$, From Lemma 2.2 in \cite[p.~39]{PK}, we have the inverse factorial expansion
\bee\label{e3e}
\frac{\g(1+2ps)}{\g(1\!+\! s) \g(1\!+\! s\!+\!\nu)}=\frac{A_0(h\kappa^\kappa)^{-s}}{(2\pi)^2}\bl\{\sum_{j=0}^{M-1}(-)^j c_j(\nu,p) \g(\kappa s+\vartheta-j)+\rho_M(s) \g(\kappa s+\vartheta-M)\br\}
\ee
for positive integer $M$, where 
\bee\label{eA0}
A_0=2\pi\kappa^{\fr-\vartheta}(2p)^\fr.
\ee
The remainder function $\rho_M(s)$ is analytic in $s$ except at the points $s=-(1+r)/(2p)$, $r=0, 1, 2 \ldots$ and is such that $\rho_M(s)=O(1)$ as $|s|\to\infty$ in $|\arg\,s|<\pi$. The coefficients $c_j(\nu,p)$ can be obtained by means of the algorithm presented in \cite[\S 2.2.4]{PK}. We have the values for $j\leq 2$
\[c_0(\nu,p)=1,\quad c_1(\nu,p)=\frac{2p-1}{12p}\{2p(1+3\nu+3\nu^2)-1\},\]
\[c_2(\nu,p)=\frac{2p-1}{288p^2}\{-1+p(-18-12\nu+12\nu^2)+p^2(36+72\nu-12\nu^2-72\nu^4-36\nu^4)\]
\bee\label{e33a}
+p^3(8+96\nu+264\nu^2+240\nu^3+72\nu^4)\}.
\ee

Then, employing the finite Fourier series
\bee\label{e33b}
\frac{\sin \pi ps}{\sin \pi s}=\sum_{r=0}^{p-1} e^{-i\kappa s\phi_r}, \qquad\phi_r:=\frac{1}{2}\pi-\frac{2\pi r}{\kappa}\quad(0\leq r\leq p-1),
\ee
we obtain
\bee\label{e33}
R_k(a;N_k)=-\frac{A_0}{8\pi}\sum_{r=0}^{p-1}\bl\{\sum_{j=0}^{M-1}(-)^j c_j(\nu,p)\,\frac{1}{2\pi i}\int_{{\cal L}_N} \frac{\g(\kappa s\!+\!\vartheta\!-\!j)}{\sin \pi(s+\nu)} (X_k e^{i\phi_r})^{-\kappa s}\,ds+{\cal R}_{r,M,N_k}\br\},
\ee
where
\[
{\cal R}_{r,M,N_k}=\frac{1}{2\pi i}\int_{{\cal L}_N}\rho_M(s)\frac{\g(\kappa s\!+\!\vartheta\!-\!M)}{\sin \pi(s+\nu)} (X_ke^{i\phi_r})^{-\kappa s}\,ds.
\] 

The remainder ${\cal R}_{r,M,N_k}$ can be estimated by means of an extended version of Lemma 2.9 in \cite[p.~75]{PK}. This extension is carried out in Appendix B, where it is shown that 
\[{\cal R}_{r,M,N_k}=O(X_k^{\vartheta-M}e^{-X_k}) \qquad (X_k\to+\infty,\ N_k\sim X_k).\] 
With the change of variable $s\to s+N_k-\nu$ in (\ref{e33}), we then find that
\[R_k(a;N_k)=-\frac{A_0 (-)^{N_k}}{8\pi}\sum_{r=0}^{p-1}\bl\{\sum_{j=0}^{M-1}(-)^j c_j(\nu,p)\,\frac{(X_ke^{i\phi_r})^\vartheta}{2\pi i}\int_{-c-\infty i}^{-c+\infty i} \frac{\g(\kappa s\!+\!\mu_k\!-\!j)}{\sin \pi s}\,(X_ke^{i\phi_r})^{-\kappa s-\mu_k}ds\]
\bee\label{e33mu}
\hspace{5cm}+O(X_k^{\vartheta-M}e^{-X_k})\br\}, \qquad \mu_k:=\vartheta+\kappa N_k-\kappa\nu,\ee
where $0<c<1$.

We now introduce the generalised {\it terminant function} $T_\om(\kappa;z)$ defined by
\bee\label{e34a}
-2e^zT_\om(\kappa;z)=\frac{1}{2\pi i}\int_{-c-\infty i}^{-c+\infty i} \frac{\g(\kappa s+\om)}{\sin \pi s}\,z^{-\kappa s-\om}ds\qquad (|\arg\,z|<p\pi/\kappa),
\ee
where $0<c<1$ and $\om\neq 0, -1, -2, \ldots\ $. Such a terminant function has been introduced previously in \cite[\S 2]{P92}. When $\kappa=1$ ($p=\f{3}{2}$) we have $T_\om(1;z)=T_\om(z)$, where $T_\om(z)$ is the standard terminant function given as a multiple of the incomplete gamma function $\g(a,z)$ by\footnote{In \cite[(2.11.11)]{DLMF} this function is denoted by $F_\om(z)$ and is expressed as a multiple of the exponential integral $E_\om(z)=z^{\om-1}\g(1-\om,z)$.}
\bee\label{e34d}
{T}_\om(z):=\frac{\g(\om)}{2\pi}\, \g(1-\om,z)=-\frac{e^{-z}}{4\pi i}\int_{-c-\infty i}^{-c+\infty i} \frac{\g(s+\om)}{\sin \pi s}\,z^{-s-\om}ds\quad (|\arg\,z|<\f{3}{2}\pi)
\ee
by (\ref{e34a}). If we define the finite sum
\bee\label{e34c}
S_j(M;X_ke^{i\phi_r}):=\sum_{j=0}^{M-1}\frac{(-)^jc_j(\nu,p)}{(X_ke^{i\phi_r})^{j-\vartheta}}\,e^{X_ke^{i\phi_r}} T_{\mu_k-j}(\kappa;X_ke^{i\phi_r}),
\ee
then the remainder can be expressed in terms of generalised terminant functions in the form
\[
R_k(a;N_k)=(-)^{N_k}\frac{A_0}{4\pi}\sum_{r=0}^{p-1}S_j(M;X_ke^{i\phi_r})+O(X_k^{\vartheta-M}e^{-X_k}).
\]

When $p$ is an even or odd integer the values of $\phi_r$ belong to the sets $\{\pm\fs\pi, \pm(\fs\pi-2\pi/\kappa), \ldots , \pm\pi/\kappa\}$ and $\{\pm\fs\pi, \pm(\fs\pi-2\pi/\kappa), \ldots , \pm 2\pi/\kappa, 0\}$, respectively.
Since we are considering $a>0$ and real $\nu$ the above expression for $R_k(a;N_k)$ can be written in the following form:
\begin{theorem}$\!\!\!.$\  For integer $p\geq 2$ we have 
\bee\label{e35a}
R_k(a;N_k)=\left\{\begin{array}{ll}{\displaystyle (-)^{N_k}\frac{A_0}{2\pi}\Re\sum_{r=0}^{(p-2)/2} S_j(M;X_ke^{i\phi_r})+O(X_k^{\vartheta-M}e^{-X_k})} & (p\ \mbox{even}) \\
\\
{\displaystyle (-)^{N_k}\frac{A_0}{4\pi}\bl\{S_j(M;X_k)+2\Re\sum_{r=0}^{(p-3)/2}  S_j(M;X_ke^{i\phi_r})\br\}+O(X_k^{\vartheta-M}e^{-X_k})} & (p\ \mbox{odd}),\end{array}\right.
\ee
where $A_0$ and $S_j(M;X_ke^{i\phi_r})$ are defined in (\ref{eA0}) and (\ref{e34c}).
\end{theorem}
 
It now remains to exploit the asymptotics of the generalised terminant function appearing in $S_j(M;Xe^{i\phi_r})$ as $X_k\to+\infty$ ($a\to0+$), where $\mu_k \sim X_k$ by (\ref{e32}). The details of the leading asymptotic behaviour of $T_\om(\kappa;z)$ are given in Appendix A.

\vspace{0.6cm}

\begin{center}
{\bf 4. \  The asymptotic form of the remainder $R_k(a;N_k)$ as $a\to0+$}
\end{center}
\setcounter{section}{4}
\setcounter{equation}{0}
\renewcommand{\theequation}{\arabic{section}.\arabic{equation}}
We now consider some specific cases and determine the leading exponential order  of the remainder function $R_k(a;N_k)$ as the parameter $a\to0+$.
\bigskip

\noindent{\bf 4.1\ \ The case $p=2$}.\ \ The associated parameters are (see (\ref{e30}) and (\ref{e33mu}))
\[\kappa=2,\qquad \mu_k=\vartheta+2N_k-2\nu,\qquad\vartheta=-\nu,\qquad \frac{A_0}{2\pi}=2^{\nu+3/2}\]
and, from (\ref{e34c}) and (\ref{e35a}), we find
\[R_k(a;N_k)=2^{\nu+3/2}(-)^{N_k}\,\Re \sum_{j=0}^{M-1}\frac{(-)^j c_j(\nu,2)}{(X_ke^{\fr\pi i})^{j-\vartheta}}\,e^{iX_k} T_{\mu_k-j}(2;X_ke^{\fr\pi i})+O(X_k^{\vartheta-M}e^{-X_k}).\]
From (\ref{e43}), we have
\[e^{iX_k} T_{\mu_k-j}(2;X_ke^{\fr\pi i})=\fs\{e^{\fr\pi i(\mu_k-j)} e^{-X_k} T_{\mu_k-j}(X_ke^{\pi i})+e^{-\fr\pi(\mu_k-j)} e^{X_k} T_{\mu_k-j}(X_k)\}\]
and hence\footnote{Use of the connection formula \cite[(6.2.45)]{PK}
$e^{\pi i\om}T_\om(ze^{\pi i})-e^{-\pi i\om}T_\om(ze^{-\pi i})=i$ shows that the form of expression given in \cite[(3.10)]{P18} agrees with (\ref{e51}).}
\[R_k(a;N_k)=2^{\nu+\fr}X_k^\vartheta e^{-X_k}\bl\{\sum_{j=0}^{M-1}\frac{(-)^j c_j(\nu,2)}{X_k^j}
\bl(\cos \pi\nu\,e^{2X_k} T_{\mu_k-j}(X_k)\hspace{5cm}\]
\bee\label{e51}
\hspace{6cm}+(-)^j\,\Re [e^{-2\pi i\nu} T_{\mu_k-j}(X_ke^{\pi i})]\br)+O(X_k^{-M})\br\}.
\ee
The coefficients $c_j(\nu,2)$ are presented in \cite{P18} for $0\leq j\leq 4$.

Application of (\ref{e44}) then shows that as $X_k\to\infty$ with $\mu_k\sim X_k$
\[T_{\mu_k-j}(X_k)=\frac{e^{-2X_k}}{2\sqrt{2\pi X_k}}\{1+O(X_k^{-1})\}\]
and
\[T_{\mu_k-j}(X_k e^{\pi i})=(-)^je^{-\pi i\mu_k}\bl(\frac{i}{2}+\frac{B_0}{\sqrt{2\pi X_k}}\{1+O(X_k^{-1})\}\br),\]
where $B_0=\f{7}{6}\!+\!\alpha_k\!-\!j$ with $\alpha_k$ specified in (\ref{e32}); see \cite[(3.17)]{O91}, \cite[(4.6)]{P18} for details. Thus we obtain
\bee\label{e51b}
R_k(a;N_k)=2^{\nu-\fr}X_k^{-\nu} e^{-X_k} \sum_{j=0}^{M-1} (-)^j\frac{ c_j(\nu,2)}{X_k^j}\bl\{\frac{\cos \pi\nu}{\sqrt{2\pi X_k}}\{1+2B_0+O(X_k^{-1})\}-\sin \pi\nu\br\}
\ee
as $a\to0+$,
where $X_k=\pi^2k^2/a$ by (\ref{e30}). A more detailed treatment of the $p=2$ case employing the asymptotic expansions of $T_{\mu_k-j}(X_k)$ and $T_{\mu_k-j}(X_ke^{\pi i})$ is given in \cite[Theorem 1]{P18}.
\bigskip

\noindent{\bf 4.2\ \ The case $p=3$}.\ \ The associated parameters are 
\[\kappa=4,\qquad \mu_k=\vartheta+4N_k-4\nu,\qquad\vartheta=-\nu,\qquad\frac{A_0}{2\pi}=2^{2\nu+1/2}\cdot 3^{1/2}\]
and, from (\ref{e34c}) and (\ref{e35a}), we find
\[R_k(a;N_k)=2^{2\nu+\fr}\cdot 3^\fr (-)^{N_k}\sum_{j=0}^{M-1}(-)^j\frac{c_j(\nu,3)}{X_k^{j-\vartheta}}\bl\{e^{X_k} T_{\mu_k-j}(4;X_k)\hspace{4cm}\]
\[\hspace{5cm}+2\Re [e^{iX_k+\fr\pi i(\vartheta-j)} T_{\mu_k-j}(4;X_ke^{\fr\pi i})]\br\}+O(X_k^{\vartheta-M}e^{-X_k}),\]
where the first few coefficients $c_j(\nu,4)$ can be obtained from (\ref{e33a}).

From (\ref{e45}) as $X_k\to\infty$ with $\mu_k\sim X_k$
\bee\label{e51a}
e^{X_k} T_{\mu_k-j}(\kappa;X_k)=\frac{e^{-X_k}}{2\sqrt{2\pi X_k}}\,\{1+O(X_k^{-1})\}
\ee
for general $\kappa\geq 1$
and, from (\ref{e46}) with $\psi_0=\f{3}{4}\pi$,
\[e^{iX_k} T_{\mu_k-j}(4;X_ke^{\fr\pi i})=\frac{e^{-X_k-\fr\pi i(\mu_k-j)}}{4\sqrt{2\pi X_k}}\,h_1(\fs\pi)\{1+O(X_k^{-1})\}\hspace{4cm}\]
\[\hspace{4cm}+\frac{i}{4} e^{-X_ke^{\pi i/4}-\frac{1}{4}\pi i(\mu_k-j)}\bl\{1-\frac{1}{2} \mbox{erfc}\,[c(\f{5}{4}\pi)(\fs X_k)^{1/2}]\br\}+O\bl(\frac{e^{-X_k}}{\sqrt{2\pi X_k}}\br),\]
where $\mbox{erfc}$ denotes the complementary error function. The quantity $c(\f{5}{4}\pi)$ is situated in the first quadrant\footnote{The locus of $c(\theta)$ as a function of $\theta$ shown in Fig.~6.3 of \cite[p.~262]{PK} should be $-c(\theta)$.} so that the argument of the above complementary error function is less than $\fs\pi$. From the asymptotic behaviour of $\mbox{erfc}\,z$  \cite[(7.12.1)]{DLMF}
\[\mbox{erfc}\,z\sim \frac{e^{-z^2}}{\sqrt{\pi}\,z}\qquad (z\to\infty\ \mbox{in}\ |\arg\,z|<\f{3}{4}\pi),\]
we therefore find from (\ref{e44e})
\[\fs \mbox{erfc}\,[c(\f{5}{4}\pi) (\fs X_k)^{1/2}]\sim \frac{e^{-\fr X_k c^2(\frac{5}{4}\pi)}}{\sqrt{2\pi X_k}\,c(\frac{5}{4}\pi)}=O\bl(\frac{e^{-X_k+X_k e^{\pi i/4}}}{\sqrt{2\pi X_k}}\br)\qquad (X_k\to\infty).
\]
Then it follows that
\[e^{iX_k} T_{\mu_k-j}(4;X_k e^{\fr\pi i})=\frac{1}{4}(-)^{N_k} e^{-X_k e^{\pi i/4}}\,e^{\pi i(\frac{1}{2}+\frac{5}{4}\nu+\frac{1}{4}j)}+O\bl(\frac{e^{-X_k}}{\sqrt{2\pi X_k}}\br).\]

Collecting together these results, employing (\ref{e51a}) with $\kappa=4$ and noting that for $M\geq 1$ the order term $O(X_k^{-M}e^{-X_k})$ can be included in the order term $O(e^{-X_k}/\sqrt{2\pi X_k})$, we finally obtain
\bee\label{e52} 
R_k(a;N_k)=2^{2\nu-\fr} 3^\fr X_k^{-\nu} \bl\{e^{-X_k/\!\surd 2} \sum_{j=0}^{M-1} (-)^j \frac{c_j(\nu,3)}{X_k^j} \sin\,\bl[\frac{X_k}{\surd 2}-\frac{3}{4}\pi\nu+\frac{1}{4}\pi j\br]+O\bl(\frac{e^{-X_k}}{\sqrt{2\pi X_k}}\br)\br\}
\ee
as $a\to0+$, where $X_k=4(2/a)^{1/2} (\pi k/3)^{3/2}$.
\bigskip

\noindent{\bf 4.3\ \ The case $p=4$}.\ \ The associated parameters are
\[\kappa=6,\qquad \mu_k=\vartheta+6N_k-6\nu,\qquad\vartheta=-\nu,\qquad\frac{A_0}{2\pi}=6^{\nu+1/2}\cdot 2^{3/2}\]
and, from (\ref{e34c}) and (\ref{e35a}), we find
\[R_k(a;N_k)= 6^{\nu+\fr} 2^\frac{3}{2}(-)^{N_k}\Re \sum_{j=0}^{M-1}(-)^j\frac{ c_j(\nu,4)}{X_k^{j-\vartheta}}\bl\{e^{\frac{1}{2}\pi i(\vartheta-j)} e^{iX_k} T_{\mu_k-j}(6;X_ke^{\fr\pi i})\hspace{2cm}\]\[\hspace{8cm}+
e^{\frac{1}{6}\pi i(\vartheta-j)} e^{X_k e^{\pi i/6}} T_{\mu_k-j}(6;X_k e^{\frac{1}{6}\pi i})\br\}.\]
From (\ref{e46}) we obtain
\[e^{X_ke^{\pi i/6}} T_{\mu_k-j}(6;X_ke^{\pi i/6})=\frac{1}{6} e^{\frac{5}{6}\pi i(\mu_k-j)} e^{-X_k}\bl\{\frac{1}{2}+\frac{1}{2} \mbox{erf}\,[c(\pi) (\fs X_k)^{1/2}\br\}+O\bl(\frac{e^{-X_k}}{\sqrt{2\pi X_k}}\br)\]\[=O(e^{-X_k})\hspace{3.8cm}\]
since $c(\pi)=0$, and from (\ref{e47})
\begin{eqnarray*}
e^{iX_k} T_{\mu_k-j}(6;iX_k)&=&\frac{i}{6}e^{-X_ke^{\pi i/3}-\frac{1}{6}\pi i(\mu_k-j)}\bl\{1-\frac{1}{2}\mbox{erfc}\,[c(\f{4}{3}\pi)(\fs X_k)^{1/2}]\br\}\\
&&+e^{-X_k-\frac{1}{2}\pi i(\mu_k-j)}\bl\{\frac{1}{2}+\frac{1}{2}\mbox{erf}\,[c(\pi)(\fs X_k)^{1/2}]\br\}+O\bl(\frac{e^{-X_k}}{\sqrt{2\pi X_k}}\br)\\
&=&\frac{i}{6}e^{-X_ke^{\pi i/3}-\frac{1}{6}\pi i(\mu_k-j)}+O(e^{-X_k}),
\end{eqnarray*}
since by the same argument employed in \S 5.2 we have $e^{-X_ke^{\pi i/3}}\mbox{erfc}\,[c(\f{4}{3}\pi)(\fs X_k)^{1/2}]=O(e^{-X_k}/\sqrt{2\pi X_k})$.

Hence we obtain
\bee\label{e53}
R_k(a;N_k)=\frac{2}{\surd 3} \br(\frac{X_k}{6}\br)^{-\nu}\bl\{e^{-X_k/2} \sum_{j=0}^{M-1}(-)^j \frac{c_j(\nu,6)}{X_k^j} \sin \bl[\frac{\surd 3}{2}X_k-\frac{2}{3}\pi\nu+\frac{1}{3}\pi j\br]+O(e^{-X_k})\br\}
\ee
as $a\to0+$, where $X_k=6(2/a)^{1/3}(\pi k/4)^{4/3}$.
\bigskip

\noindent{\bf 4.4\ \ The case $p=5$}.\ \ In our final example the associated parameters are
\[\kappa=8,\qquad \mu_k=\vartheta+8N_k-8\nu,\qquad\vartheta=-\nu,\qquad\frac{A_0}{2\pi}=2^{3\nu+2}\cdot 5^{1/2}\]
and from (\ref{e34c}) and (\ref{e35a}), we find
\[R_k(a;N_k)= 2^{3\nu+1} \surd 5 (-)^{N_k}\sum_{j=0}^{M-1} (-)^j \frac{c_j(\nu,5)}{X_k^j}\bl\{e^{X_k}T_{\mu_k-j}(8;X_k)\hspace{4cm}\]
\bee\label{e54}
\hspace{4cm}+2\Re \sum_{r=0}^1 e^{X_ke^{i\phi_r+i\phi_r(\vartheta-j)}}T_{\mu_k-j}(8;X_ke^{i\phi_r})\br\}
+O(X_k^{\vartheta-M}e^{-X_k}),
\ee
where $\phi_0=\fs\pi$ and $\phi_1=\f{1}{4}\pi$.

From (\ref{e51a}) with $\kappa=8$, the first generalised terminant in the above expression is $O(e^{-X_k}/\sqrt{2\pi X_k})$. From (\ref{e46}) and (\ref{e47}) with $\psi_0=\f{7}{8}\pi$, $\psi_1=\f{5}{8}\pi$ we find
\[e^{X_ke^{\pi i/4}}T_{\mu_k-j}(8;X_ke^{\frac{1}{4}\pi i})=\frac{i}{8} e^{-X_ke^{\pi i/8}-\frac{1}{8}\pi i(\mu_k-j)}\bl\{1-\frac{1}{2}\mbox{erfc}\,[c(\f{9}{8}\pi) (\fs X_k)^{1/2}]\br\}+O\bl(\frac{e^{-X_k}}{\sqrt{2\pi X_k}}\br)\]
\[=\frac{i}{8} e^{-X_ke^{\pi i/8}-\frac{1}{8}\pi i(\mu_k-j)}+O\bl(\frac{e^{-X_k}}{\sqrt{2\pi X_k}}\br)\]
and
\begin{eqnarray*}
e^{iX_k} T_{\mu_k-j}(8;X_ke^{\fr\pi i})&=&\frac{i}{8} e^{-X_ke^{3\pi i/8}-\frac{1}{8}\pi i(\mu_k-j)}\bl\{1-\frac{1}{2}\mbox{erfc}\,[c(\f{11}{8}\pi) (\fs X_k)^{1/2}]\br\}\\
&&\!\!\!\!+\frac{i}{8} e^{-X_ke^{\pi i/8}-\frac{3}{8}\pi i(\mu_k-j)}\bl\{1-\frac{1}{2}\mbox{erfc}\,[c(\f{9}{8}\pi) (\fs X_k)^{1/2}]\br\}+O\bl(\frac{e^{-X_k}}{\sqrt{2\pi X_k}}\br)\\
&=&\frac{i}{8}\bl\{e^{-X_ke^{3\pi i/8}-\frac{1}{8}\pi i(\mu_k-j)}+e^{-X_ke^{\pi i/8}-\frac{3}{8}\pi i(\mu_k-j)}\br\}+O\bl(\frac{e^{-X_k}}{\sqrt{2\pi X_k}}\br).
\end{eqnarray*}
Substitution of these leading asymptotic forms in (\ref{e54}) then yields after some routine algebra
\[R_k(a;N_k)=\frac{\surd 5}{2}\bl(\frac{X_k}{8}\br)^{\!-\nu}\sum_{j=0}^{M-1}(-)^j\frac{c_j(\nu,8)}{X_k^j}
\bl\{e^{-X_k \cos 3\pi/8} \sin\bl[X_k \sin \f{3}{8}\pi-\frac{5}{8}\pi\nu+\frac{3}{8}\pi j\br]\hspace{2cm}\]
\bee\label{e55}
\hspace{3cm}+2\cos \pi\nu\,e^{-X_k\cos \pi/8} \sin\bl[X_k\sin \f{1}{8}\pi-\frac{15}{8}\pi\nu+\frac{1}{8}\pi j\br]\br\}
+O\bl(\frac{e^{-X_k}}{\sqrt{2\pi X_k}}\br)
\ee
as $a\to0+$, where $X_k=8(2/a)^{1/4} (\pi k/5)^{5/4}$.
\vspace{0.6cm}

\begin{center}
{\bf 5. \  Concluding remarks}
\end{center}
\setcounter{section}{5}
\setcounter{equation}{0}
\renewcommand{\theequation}{\arabic{section}.\arabic{equation}}
When $p=2$ the optimally truncated remainder $R_k(a;N_k)$ defined in (\ref{e35a}) is exponentially small of $O(e^{-X_k})$ as $a\to0+$. This case is treated more fully in \cite{P18}. When $p=3, 4$ a second subdominant exponential contribution, in addition to the $O(e^{-X_k})$ contribution, appears given by $O(e^{-X_k/\surd 2})$ and $O(e^{-X_k/2})$, respectively. These additional contributions, although subdominant, are less recessive than $e^{-X_k}$. When $p=5$, there are two additional contributions of $O(e^{-X_k \cos 3\pi/8})$ and 
$O(e^{-X_k \cos \pi/8})$, which together with the $O(e^{-X_k})$
contribution form a sequence of increasingly subdominant exponential contributions.

Examination of (\ref{e35a}) reveals that this pattern continues with progressively more exponential contributions of decreasing subdominance appearing as $p$ increases. It can be shown that the subdominant contributions in the remainder $R_k(a;N_k)$ at optimal truncation are controlled by
\[\exp \,\bl[-X_k \cos \bl(\frac{\pi(p-2\ell)}{\kappa}\br)\br]\qquad (\ell=1, 2, \ldots , \lceil\fs p\rceil-1).\]
Consequently, as $p$ is allowed to increase the number of exponentially small contributions increases by one each time $p$ equals an odd integer. A similar observation was made in the discussion of the Euler-Jacobi series 
$\sum_{n\geq 1} e^{-an^p}$ \cite{Kow}, \cite[\S 8.1.2]{PK}. In \cite{PK}, the appearance of additional exponentially small contributions was seen to be associated with the appearance of additional saddle points in the asymptotic treatment of the series by
the saddle-point method applied to a Laplace-type integral representation.
\vspace{0.6cm}

\vspace{0.6cm}

\begin{center}
{\bf Appendix A: The asymptotics of the generalised terminant function $T_\om(\kappa;z)$}
\end{center}
\setcounter{section}{1}
\setcounter{equation}{0}
\renewcommand{\theequation}{\Alph{section}.\arabic{equation}}
The generalised terminant function $T_\om(\kappa;z)$ is defined by the integral
\bee\label{e41}
-2z^\om e^z T_\om(\kappa;z)=\frac{1}{2\pi i}\int_{-c-\infty i}^{-c+\infty i}\frac{\g(\kappa s+\om)}{\sin \pi s}\,z^{-\kappa s}ds\qquad(|\arg\,z|<p\pi/\kappa),
\ee
where we recall that $\kappa=2(p-1)$, $\om\neq 0, -1, -2, \ldots$ and $0<c<1$. For integer $p\geq 2$, we have from (\ref{e33b})
\bee\label{e42}
\frac{\sin \pi\kappa s}{\sin \pi s}=\sum_{r=0}^{\kappa-1} e^{-i\kappa s\psi_r}, \qquad\psi_r:=\pi-\frac{(2r+1)\pi }{\kappa}\quad(0\leq r\leq \kappa-1).
\ee 
Then, with $c$ chosen to satisfy $0<c<1/\kappa$, the right-hand side of (\ref{e41}) can be expressed as
\begin{eqnarray*}
&&\sum_{r=0}^{\kappa-1} \frac{1}{2\pi i}\int_{-c-\infty i}^{-c+\infty i} \frac{\g(\kappa s+\om)}{\sin \pi\kappa s}\,(ze^{i\psi_r})^{-\kappa s} ds\\
&=&\frac{1}{\kappa}\sum_{r=0}^{\kappa-1} \frac{1}{2\pi i}\int_{-c'-\infty i}^{-c'+\infty i} \frac{\g(s+\om)}{\sin \pi s}\,(ze^{i\psi_r})^{-s} ds\qquad(0<c'<1)\\
&=&\frac{1}{\kappa}\sum_{r=0}^{\kappa-1}(ze^{i\psi_r})^\om e^{ze^{i\psi_r}} T_\om(ze^{i\psi_r})
\end{eqnarray*}
by (\ref{e34d}). Hence we obtain \cite[\S 2]{P92}
\bee\label{e43}
e^z T_\om(\kappa;z)=\frac{1}{\kappa}\sum_{r=0}^{\kappa-1}e^{i\om\psi_r} e^{ze^{i\psi_r}} T_\om(ze^{i\psi_r}),
\ee
when $\kappa=2(p-1)$ and integer $p\geq 2$. We remark that the values of $\psi_r$ belong to the set $\{\pm(\pi-\pi/\kappa), \pm(\pi-3\pi/\kappa), \ldots , \pm\pi/\kappa\}$.

The asymptotic expansion of the terminant function $T_\om(z)$ for large $\om$ and complex $z$, when $\om\sim |z|$, has been discussed in detail by Olver in \cite{O91}; see also \cite[Section 2.11]{DLMF} and the detailed account in \cite[pp.~259--265]{PK}. It is found that\footnote{It has been shown in \cite[\S 4]{O91} that the second expansion holds in the wider sector $(-\pi,3\pi)$.}
\bee\label{e44}
T_\om(z)=\left\{\begin{array}{ll}
\displaystyle{\frac{e^{-i\om\theta}}{1+e^{-i\theta}}\,\frac{e^{-z-|z|}}{\sqrt{2\pi |z|}}\{1+O(z^{-1})\}} & |\arg\,z|\leq \pi-\delta, \delta>0\\
\displaystyle{ie^{-\pi i\om}\bl\{\frac{1}{2}+\frac{1}{2} \mbox{erf}\,[c(\theta) (\fs|z|)^{1/2}]\br\}+O\bl(\frac{e^{-z-|z|}}{\sqrt{2\pi|z|}}\br)} & \delta\leq\arg\,z\leq2\pi-\delta \end{array}\right.
\ee
valid as $|z|\to\infty$ when $\om\sim |z|$, where $\theta=\arg\,z$. A conjugate exansion holds for $-2\pi+\delta\leq\arg\,z\leq-\delta$.
The quantity $c(\theta)$ is defined by
\bee\label{e44e}
\fs c^2(\theta)=1+i(\theta-\pi)-e^{i(\theta-\pi)}
\ee
and has the expansion 
$c(\theta)=\theta-\pi+\f{1}{6}i(\theta-\pi)^2+O((\theta-\pi)^3)$ in the neighbourhood of $\theta=\pi$. Then from (\ref{e43}) and (\ref{e44}) we can construct the leading asymptotics of $T\om(\kappa;z)$ when $\om\sim|z|$.

Let us define the quantity
\[h_m(\theta)=\sum_{r=m}^{\kappa-1}\frac{1}{1-\lambda_r e^{-i\theta}},\quad\lambda_r=e^{(2r+1)\pi i/\kappa}\qquad(0\leq m\leq\kappa-1),\]
where $h_0(\theta)=\kappa/(1+e^{-i\kappa\theta})$. Then, when $|\theta|<\pi/\kappa$, it is seen that $|\arg\,(ze^{i\psi_r})|<\pi$ so that the first expression in (\ref{e44}) applies in the sum of terminant functions on the right-hand side of (\ref{e43}). Hence we find \cite[(2.10)]{P92}
\begin{eqnarray}
e^z T_\om(\kappa;z)&=&\frac{e^{-|z|}}{\kappa\sqrt{2\pi |z|}} \sum_{r=0}^{\kappa-1}\frac{e^{-i\om\theta}}{1-\lambda_r e^{-i\theta}}\,\{1+O(|z|^{-1})\}\nonumber\\
&=&\frac{e^{-|z|}}{\sqrt{2\pi |z|}}\,\frac{e^{-i\om\theta}}{1+e^{-i\kappa\theta}}\,\{1+O(|z|^{-1})\}\qquad(|\theta|<\pi/\kappa)\label{e45}
\end{eqnarray}
as $|z|\to\infty$ when $\om\sim|z|$.

When $\pi/\kappa\leq\theta<3\pi/\kappa$ we have $|\arg\,(ze^{i\psi_r})|<\pi$ for $1\leq r\leq \kappa-1$ and $\pi\leq\arg\,(ze^{i\psi_0})<\pi+2\pi/\kappa$. Hence
\begin{eqnarray}
e^z T_\om(\kappa;z)&=&\frac{1}{\kappa}e^{ze^{i\psi_0}+i\om\psi_0} T_\om(ze^{i\psi_0})+\frac{e^{-|z|-i\om\theta}}{\kappa\sqrt{2\pi|z|}}\,h_1(\theta)\{1+O(|z|^{-1})\}\nonumber\\
&=&\frac{i}{\kappa}e^{-ze^{-\pi i/\kappa}-\pi i\om/\kappa}\bl\{\frac{1}{2}+\frac{1}{2}\,\mbox{erf}\,[c(\theta+\psi_0) (\fs|z|)^{1/2}]\br\}+O\bl(\frac{e^{-|z|}}{\sqrt{2\pi|z|}}\br)\nonumber\\
&&+\frac{e^{-|z|-i\om\theta}}{\kappa\sqrt{2\pi|z|}}\,h_1(\theta)\{1+O(|z|^{-1})\},\label{e46}
\end{eqnarray}
where $c(\theta+\psi_0)\simeq \theta-\pi/\kappa$ near $\theta=\pi/\kappa$. A conjugate relation holds when $-3\pi/\kappa<\theta\leq-\pi/\kappa$.

When $3\pi/\kappa\leq\theta<5\pi/\kappa$ we have $|\arg\,(ze^{i\psi_r})|<\pi$ for $2\leq r\leq \kappa-1$ and
$\pi+2\pi/\kappa\leq\arg\,(ze^{i\psi_0})<\pi+4\pi/\kappa$ and $\pi\leq\arg\,(ze^{i\psi_1})<\pi+2\pi/\kappa$. Thus we have the result
\begin{eqnarray}
e^zT_\om(\kappa;z)&=&\frac{i}{\kappa}e^{-ze^{-\pi i/\kappa}-\pi i\om/\kappa}\bl\{\frac{1}{2}+\frac{1}{2}\,\mbox{erf}\,[c(\theta+\psi_0) (\fs|z|)^{1/2}]\br\}\nonumber\\
&&+\frac{i}{\kappa}e^{-ze^{-3\pi i/\kappa}-3\pi i\om/\kappa}\bl\{\frac{1}{2}+\frac{1}{2}\,\mbox{erf}\,[c(\theta+\psi_1) (\fs|z|)^{1/2}]\br\}+O\bl(\frac{e^{-|z|}}{\sqrt{2\pi|z|}}\br)\nonumber\\
&&+\frac{e^{-|z|-i\om\theta}}{\kappa\sqrt{2\pi|z|}}\,h_2(\theta)\{1+O(|z|^{-1})\},\label{e47}
\end{eqnarray}
where $c(\theta+\psi_1)\simeq \theta-3\pi/\kappa$ near $\theta=3\pi/\kappa$. A conjugate relation holds when $-5\pi/\kappa<\theta\leq-3\pi/\kappa$.

Other ranges of $\theta=\arg\,z$ that correspond to increasingly more terminant functions in (\ref{e43}) possessing an argument greater (in modulus) than $\pi$, provided $p$ is such that the range does not exceed the domain of validity in (\ref{e44}).  
\vspace{0.6cm}

\begin{center}
{\bf Appendix B: A bound on the remainder ${\cal R}_{r,M,N_k}$ }
\end{center}
\setcounter{section}{2}
\setcounter{equation}{0}
\renewcommand{\theequation}{\Alph{section}.\arabic{equation}}
In \cite[Lemma 2.9]{PK} the order of the integral
\[J(z)=\frac{1}{2\pi i}\int_{-c+n-\infty i}^{-c+n+\infty i} \rho_M(s) \g(s+\alpha)\,\frac{z^{-s}}{\sin \pi s}\qquad (|\arg\,z|<\f{3}{2}\pi,\ 0<c<1)\] 
is considered when the integer $n$ is chosen to satisfy $n\sim|z|$ for $|z|\to\infty$. The quantity $\rho_M(s)=O(1)$ as $|s|\to\infty$ in $|\arg\,s|<\pi$ and the parameter $\alpha$ is a bounded arbitrary real constant. It was established that as $|z|\to\infty$
\[|J(z)|=\left\{\begin{array}{ll}O(z^{\alpha-\fr}e^{-|z|}) & (|\arg\,z|<\pi)\\
\\
O(z^\alpha e^{-|z|}) & (|\arg\,z|\leq\pi).\end{array}\right.\]
Here we extend the order estimate to the range $\pi\leq|\arg\,z|<\f{3}{2}\pi$ to yield the following lemma:
\vspace{0.2cm}

\noindent
{\bf Lemma 1.}\ \ \ {\it For $n\sim |z|$ and $|z|\to\infty$ we have the order estimate
\bee\label{a1}
|J(z)|=O(z^\alpha e^{-|z|}(\sec \varphi)^n) \qquad (\pi\leq|\arg\,z|<\f{3}{2}\pi),
\ee
where $\varphi=|\arg\,z|-\pi$.}
\vspace{0.2cm}

\noindent
{\bf Proof.}\ \ Let $\theta=\arg\,z$ and $\varphi=\theta-\pi$; then for $\pi\leq\theta<\f{3}{2}\pi$ we have $0\leq\varphi<\fs\pi$.
Set $\sigma=n+\alpha-c$ and use the bound $|\rho_M(s)|<K$ for $|s|>K'$, where $K$, $K'$ are positive constants. Then for $t\in(-\infty,\infty)$, we have 
\[\frac{e^{\theta t}}{(\cosh^2 \pi t-\cos^2 \pi c)^{1/2}}
\leq\frac{e^{\theta t}}{\cosh \pi t\,\sin \pi c}\leq\frac{2e^{\theta t-\pi|t|}}{\sin \pi c}\]
when $0<c<1$, where
\[e^{\theta t-\pi|t|}=e^{\varphi t}\ \  (t\geq 0),\quad e^{\theta t-\pi|t|}=e^{-(2\pi+\varphi)|t|}=e^{\varphi t} e^{-2\pi|t|}\ \  (t<0).\]
Hence, with $s=n-c+it$ we have
\[|J(z)|<\frac{K|z|^{-n+c}}{2\pi}\int_{-\infty}^\infty |\g(\sigma+it)|\,\frac{e^{\theta t}}{(\cosh^2 \pi t-\cos^2 \pi c)^{1/2}}\,dt\]
\[<\frac{K|z|^{-n+c}}{\pi\sin \pi c}\int_{-\infty}^\infty |\g(\sigma+it)|\,e^{\varphi t}dt.\hspace{2cm}\]

By \cite[Lemma 2.6]{PK} the last integral is $O(e^{-\sigma} \sigma^\sigma (\sec \varphi)^{\sigma+\fr})$ as $\sigma\to+\infty$. Hence with $n\sim |z|$ we obtain as $|z|\to\infty$
\[|J(z)|=O(z^\alpha e^{-|z|} (\sec \varphi)^n) \qquad (\pi\leq\arg\,z<\f{3}{2}\pi).\]
It is seen that the above estimate on the bound breaks down as $\varphi\to\fs\pi$ ($\theta\to\f{3}{2}\pi$).
A similar argument can be given when $-\f{3}{2}\pi<\arg\,z\leq -\pi$. This establishes (\ref{a1}). \hfill $\Box$
\bigskip

We can now employ Lemma 1 to estimate the remainder $R_{r,M,N_k}$ in (\ref{e33a}), where upon use of the above-mentioned bound for $\rho_M(s)$ and replacement of $s$ by $s-\nu$ we find
\[|R_{r,M,N_k}|<\frac{X_k^{\kappa\nu}}{2\pi} \left|\int_{-c+N_k-\infty i}^{-c+N_k+\infty i} \rho_M(s\!-\!\nu)\,\frac{\g(\kappa s\!+\!\alpha)}{\sin \pi s}\,(X_ke^{i\phi_r})^{-\kappa s} ds\right|\qquad \alpha:=\vartheta-\kappa\nu-M.\]
Use of the expansion in (\ref{e42}) then leads to
\[|R_{r,M,N_k}|<\frac{X_k^{\kappa\nu}}{2\pi}\sum_{j=0}^{\kappa-1} \left|\int_{-c+N_k-\infty i}^{-c+N_k+\infty i} \rho_M(s\!-\!\nu)\,\frac{\g(\kappa s\!+\!\alpha)}{\sin \pi\kappa s}\,(X_ke^{i(\phi_r+\psi_j)})^{-\kappa s} ds\right|,\]
where $c$ is chosen to satisfy $0<c<1/\kappa$ and $\psi_j=\pi-(2j+1)\pi/\kappa$. It is easily verified from (\ref{e33b}) and (\ref{e42}) that $|\phi_r+\psi_j|\leq \f{3}{2}\pi-\pi/\kappa$, whereupon from Lemma 1 with $N_k\sim X_k$ we obtain
\bee\label{a2}
|R_{r,M,N_k}|=X_k^{\kappa\nu}\, O(X_k^\alpha e^{-X_k})=O(X_k^{\vartheta-M}e^{-X_k}) \qquad (X_k\to+\infty).
\ee

\end{document}